\documentclass[11pt]{amsart}
\usepackage{amsmath}
\usepackage[all]{xy}
\usepackage{amssymb}
\usepackage{mathrsfs}
\usepackage{cite}
\usepackage{hyperref}
\usepackage{amsfonts}
\usepackage{mathdots}

\theoremstyle{plain}
\newtheorem{thm}{Theorem}[section]
\newtheorem{prop}[thm]{Proposition}
\newtheorem{lem}[thm]{Lemma}
\newtheorem{cor}[thm]{Corollary}

\theoremstyle{remark}
\newtheorem{rem}[thm]{Remark}

\theoremstyle{definition}
\newtheorem{defn}[thm]{Definition}

\newcommand{\cA}{\mathcal{A}}
\theoremstyle{conjecture}
\newtheorem{conj}[thm]{Conjecture}
\newcommand{\Ku}{\mathcal{K}u}

\def\HH{\mathrm{HH}}
\def\Hom{\mathrm{Hom}}
\def\Id{\mathrm{Id}}

\DeclareMathOperator\oh{\mathcal{O}}
\newcommand{\bv}{\mathbf{v}}

\title{Kuznetsov's Fano threefold conjecture via Hochschild-Serre algebra}

\date{}

\begin{document}

\author{Xun Lin and Shizhuo Zhang}

\begin{abstract}
  Let $Y$ be a smooth quartic double solid regarded as a degree 4 hypersurface of the weighted projective space $\mathbb{P}(1,1,1,1,2)$. We study the multiplication of the Hochschild-Serre algebra of its Kuznetsov component $\Ku(Y)$ via matrix factorization. As an application, we give a new disproof of Kuznetsov's Fano threefold conjecture. 
\end{abstract}



\address{Max Planck Institute for Mathematics, Vivatsgasse 7, 53111 Bonn, Germany}
\email{xlin@mpim-bonn.mpg.de, lin-x18@tsinghua.org.cn}

\address{Simons Laufer Mathematical Sciences Institute, 17 Gauss Way, Berkeley, CA 94720, USA}
\address{Institut de Mathématiqes de Toulouse, UMR 5219, Université de Toulouse, Université Paul Sabatier, 118 route de
Narbonne, 31062 Toulouse Cedex 9, France}
\address{MCM, Academy of Mathematics and System Science, Chinese Academy of Sciences, Beijing 100190}
\email{shizhuozhang@msri.org,shizhuo.zhang@math.univ-toulouse.fr}

\maketitle

\section{Introduction}
Let $X$ be Fano variety whose semi-orthogonal decomposition for bounded derived category is given by 
$$D^b(X)=\langle\Ku(X),E_1,\ldots,E_n\rangle,$$
where $E_1,\ldots,E_n$ is an exceptional collection of vector bundles over $X$ and $\Ku(X)$ is the right orthogonal complement of the collection, called Kuznetsov component. It has been widely believed that the Kuznetsov component encodes the essential birational geometric information of the Fano varieties. Thus extracting geometric information from Kuznetsov components is an important step to understand geometry of Fano varieties. There are numerous way to extract information from Kuznetsov components, which we briefly recall as follows. 

\subsection{Stability conditions in Kuznetsov components and moduli space theoretical approach}\label{section_moduli_space}
One of the most interesting class of Fano varieties are smooth Fano threefolds of Picard rank one of index one and two, whose deformation classes are completely classified in \cite{shafa}. In the paper \cite{bayer2017stability}, the authors construct a stability condition $\sigma$ in $\Ku(X)$ for any such Fano variety $X$. Denote by $\mathcal{N}(\Ku(X))$ the numerical Grothendieck group and fix a numerical class $\bv\in\mathcal{N}(\Ku(X))$ and consider the Bridgeland moduli space $\mathcal{M}_{\sigma}(\Ku(X),\bv)$ of (semi)stable object with respect to $\sigma$ in $\Ku(X)$ of numerical class $\bv$. The numerical character $\bv$ is appropriately chosen such that the corresponding moduli space reconstructs Fano variety of rational curves on $X$, which is used to reconstruct (birational) isomorphism class of Fano varieties(cf. \cite{bernardara2012categorical},\cite{pertusi:some-remarks-fano-threefolds-index-two}, \cite{guo2022conics}). 

\subsection{Topological K-theory of admissible subcategory and Hodge theoretical approach}\label{section_Hogde_theoretical}
Let $\cA\subset D^b(X)$ be an admissible subcategory of bounded derived category of a smooth projective variety $X$. The topological K-theory \cite{blanc_2016} of dg categories over $\mathbb{C}$ is an additive invariant 
$$\mathrm{K}^{top}_{1}: \mathrm{dg}-cat\rightarrow \mathbb{Z}-mod.$$
 with Chern character map 
 $$\mathrm{ch}^{top}: \mathrm{K}^{top}_{1}(\mathcal{A}_{dg})\rightarrow \mathrm{HP}_{1}(\mathcal{A}_{dg}).$$
 Furthermore $\mathrm{K}^{top}_{1}(D^{perf}_{dg}(X))\otimes \mathbb{C}\cong \mathrm{H}^{\text{odd}}(X,\mathbb{C})$, and $\mathrm{ch}^{top}$ is the usual Chern character. In particular, the natural splitting from that of $X$ gives a weight one Hodge structure for topological $K$-group $\mathrm{K}^{top}_{1}(\mathcal{A})_{tf}$. Namely, the topological Chern character induces,
$$\mathrm{K}^{top}_{1}(\mathcal{A})_{tf}\otimes \mathbb{C}\cong^{\mathrm{ch}^{top}} \mathrm{HP}_{1}(\mathcal{A})\cong \mathrm{HN}_{-1}(\mathcal{A})\oplus \overline{\mathrm{HN}_{-1}(\mathcal{A})}.$$
Thus, we have a complex torus associated to this weight one Hodge structure. More explicitly,
$$\mathrm{J}(\mathcal{A})=\frac{\mathrm{HP}_{1}(\mathcal{A})}{\mathrm{HN}_{-1}(\mathcal{A})+\operatorname{ch^{top}}(\mathrm{K}^{top}_{1}(\mathcal{A}))}.$$
In the case of $X$ is a smooth Fano threefold with $\cA$ being the Kuznetsov component $\Ku(X)$ as the orthogonal complement of an exceptional collection of vector bundles, by \cite[Lemma 3.9]{jacovskis2021categorical} $\mathrm{J}(\Ku(X))\cong J(X)$ as polarised abelian varieties(which was sketched in \cite[Section 5]{perry2020integral}). Similar construction is generalized for any smooth and proper dg category and even to arbitrary dg category in \cite{CMHLZZ2024categorical} and \cite{LXZ2024Hodge}. On the other hand, in the similar spirit, topological K-theory and noncommutative Hodge theory(cf. \cite{blanc_2016} and \cite{perry2020integral}) is applied to admissible subcategory of Fano fourfolds in \cite{bayer2022kuznetsov} and surfaces in \cite{dell2023cyclic} to recover Mukai lattice for K3 category and primitive cohomology for surfaces respectively. As application, (birational) categorical Torelli theorem are proved for many varieties. 

\subsection{Hochschild-Serre algebra and algebraic approach}\label{section_algebraic_approach}
The Hochschild cohomology is an algebra, and Hochschild homology is a graded module over this algebra.
We now define a bi-graded algebra that contains Hochschild cohomology and Hochschild homology and encodes the algebra structure of Hochschild cohomology and the module structure of Hochschild homology over Hochschild cohomology. Let $\cA$ be a smooth and proper dg category and $S$ be the Serre functor of $\cA$. One can naturally attach a bi-graded algebra $$\mathcal{A}_{S}=\bigoplus_{m,n\in \mathbb{Z}} \Hom(\Id,S^{m}[n])$$ with multiplication map  $$\xymatrix@C=2cm{\Hom(\Id,S^{m_1}[n_1])\times \Hom(\Id,S^{m_2}[n_2])\ar[r]^{\times}&\Hom(\Id,S^{m_1+m_2}[n_1+n_2])}$$ given by the composition
 $$\xymatrix@C=3.5cm{\Id\ar[r]^{b}&\Id\circ S^{m_2}[n_2]\ar[r]^{a\circ \Id}&S^{m_1}[n_1]\circ S^{m_2}[n_2]=S^{m_1+m_2}[n_1+n_2]
 },$$ for $(a,b)\in \Hom(\Id,S^{m_1}[n_1])\times \Hom(\Id,S^{m_2}[n_2])$. It was studied in \cite{orlov2003derived} and \cite{cualduararu2005derived}, \cite{caldararu2003mukai} independently when $\cA$ is bounded derived category $D^b(X)$ of coherent sheaves on a smooth projective variety $X$, where they prove basic property of this algebra. Moreover, in \cite{bondal2001reconstruction}, the author uses a sub-algebra of $D^b(X)_S$, which is isomorphic to anti-canonical ring of a smooth Fano variety $X$ to reconstruct the variety itself. Recently this algebra is revisited in \cite{belmans2023hochschild} and \cite{lin2023serre} under the name \emph{Hochschild-Serre algebra} for admissible subcategory of $D^b(X)$. In particular, the authors of \cite{lin2023serre} establish a sub-algebra of $\Ku(X)_S$ in the case of smooth hypersurface of degree $d$ in $\mathbb{P}^n$, which recovers the Jacobian ring of $X$ if $\mathrm{gcd}(d,n+1)=1$. Thus a categorical Torelli theorem is proved for those hypersurfaces. 

 \subsection{Kuznetsov's Fano threefold conjecture}
 Denote by $\mathcal{MF}^i_d$ the moduli space of smooth Fano threefold of index $i$ and degree $d$. In \cite[Conjecture 3.7]{kuznetsov2009derived}, the author proposed a surprising conjecture relating the non-trivial admissible subcategories of two families of smooth Fano threefolds. 

\begin{conj}\label{conjecture_Kuz_Fano_threefold}
There is a correspondence $\mathcal{Z}_d\subset\mathcal{MF}^2_d\times\mathcal{MF}^{1}_{4d+2}$, such that for any pair $(Y_d,X_{4d+2})\in\mathcal{Z}_d$, there is an equivalence of categories
$$\Ku(Y_d)\simeq\Ku(X_{4d+2}).$$
\end{conj}

The conjecture is proved for $d=3,4$ and $5$ in \cite{kuznetsov2009derived}. For the remaining cases, there was a lot of evidence suggesting that the conjecture might be false. Thus instead of proving this conjecture, people disproved it. To do this, the natural idea would be looking at the information(moduli spaces, Hodge theory, algebra etc.)extracting from $\Ku(Y_d)$ and $\Ku(X_{4d+2})$ respectively and then show that they are different. Indeed, in \cite{Zhang2020Kuzconjecture}, the author adopts the moduli theoretical approach described in Section~\ref{section_moduli_space} to study particular Bridgeland moduli spaces canonically constructed from $\Ku(Y_2)$ and $\Ku(X_{10})$ respectively and shows that they are not isomorphic to each other. Independently, in \cite{bayer2022kuznetsov}, the authors apply the Hodge theoretical approach in Section~\ref{section_Hogde_theoretical}. They look at the Mukai-Hodge lattice of two K3 categories constructed from equivariant categories of $\Ku(Y_2)$ and $\Ku(X_{10})$ respectively and show that the Hodge isometry does not exist. In the current paper, we adopt another perspective described in \ref{section_algebraic_approach}. Namely, we look at Hochschild-Serre algebra of dg-enhancement of Kuznetsov component of quartic double solid and Gushel-Mukai threefold, which are the Fano threefolds appearing in the $d=2$ case of Conjecture~\ref{conjecture_Kuz_Fano_threefold}. It turns out that they are not isomorphic to each other and the proof is very simple. 

\subsection{Main Results}
Let $Y$ be a smooth quartic double solid with the geometric involution $\iota$ and $X$ be a smooth Gushel-Mukai threefold. Denote by $\Ku(Y)_S,\Ku(X)_S$ the Hochschild-Serre algebra of dg-enhancement of $\Ku(Y)$ and $\Ku(X)$ respectively. Note that $S_{\Ku(Y)}=\iota\circ [2]$ by \cite[Section 3]{altavilla2022moduli}. Then, consider the multiplication

\begin{align}\mathrm{Hom}(\mathrm{Id},S_{\Ku(Y)}^2[-2])\times\mathrm{Hom}(\mathrm{Id},S_{\Ku(Y)}[-1])\rightarrow\mathrm{Hom}(\mathrm{Id},S_{\Ku(Y)}^3[-3])\cong\mathrm{HH}_1(\Ku(Y)),\end{align}
and associated map 
$$\gamma_{Y}:\mathrm{HH}^2(\Ku(Y))\rightarrow\mathrm{Hom}(\mathrm{HH}_{-1}(\Ku(Y)),\mathrm{HH}_1(\Ku(Y))).$$ Then we show 

\begin{thm}\label{thm_main_result}
The kernel of the map $\gamma_{Y}$ is one dimensional.
\end{thm}

On the other hand, from \cite[Theorem 4.6]{jacovskis2022infinitesimal}, we know the map 
$$\gamma_X:\mathrm{HH}^2(\Ku(X))\rightarrow\mathrm{Hom}(\mathrm{HH}_{-1}(\Ku(X)),\mathrm{HH}_1(\Ku(X)))$$ is injective for all ordinary Gushel-Mukai threefold $X$. In fact, the injectivity holds for special
Gushel–Mukai threefolds as well, as we will show in Lemma~\ref{lemma_injectivity_gamma}. Then by \cite[Theorem 4.8]{jacovskis2022infinitesimal}, $\mathrm{Ker}(\gamma_Y)\cong\mathrm{Ker}(\gamma_X)=0$, which is a contradiction. Thus we have 

\begin{cor}
For any Gushel Mukai threefold $X$ and quartic double solid $Y$, the categories $\Ku(X)$ and $\Ku(Y)$ are never equivalent. In particular, the Conjecture~\ref{conjecture_Kuz_Fano_threefold} for $d=2$ fails. 
\end{cor}

\subsection{Organization of the paper}
In Section~\ref{section_intro_matrix_factorizations}, we recall the terminology of category of graded matrix factorization $\mathrm{Inj}_{\mathrm{coh}}(\mathbb{A}^{n+1},\mathbb{C}^{\ast},\oh(d),\omega)$ with $\mathbb{C}^{\ast}$-action on $\mathbb{A}^{n+1}$ of weight $(a_0,\ldots,a_n)$. Then we describe the multiplications of Hochschild-Serre algebra for the matrix factorization. In Section~\ref{section_main_theorem_proved}, we describe the multiplication of Hochschild-Serre algebra for Kuznetsov component of a smooth quartic double solid and prove 
Theorem~\ref{thm_main_result}, as a corollary, we disproof Kuznetsov's Fano threefold conjecture.


\subsection{Acknowledgement}
We would like to thank Marcello Bernardara, Pieter Belmans, Will Donovan, Junwu Tu for useful conversation on related topics. We also would like to thank the referee for the careful reading and for providing detailed comments. SZ is supported by ANR project FanoHK, grant ANR-20-CE40-0023, Deutsche Forschungsgemeinschaft under Germany's Excellence Strategy-EXC-2047$/$1-390685813. Part of the work was finished when XL and SZ visited the Max-Planck Institute for Mathematics, Hausdorff Institute for Mathematics and Morningside Center of Mathematics, Chinese Academy of Sciences. They are grateful for excellent working conditions and hospitality.

\section{dg category of graded matrix factorizations}\label{section_intro_matrix_factorizations}
In this section, we recall the terminology of $dg$-category of matrix factorization. We follow the context in \cite{BFK}. We refer the reader to \cite{keller2006differential} for the basic of $dg$ categories. 
Denote by $\operatorname{Hqe(dg-cat)}$ the localized $\operatorname{dg-cat}$ with respect to the quasi-equivalences of $dg$ categories.   Let $(X,G,L,\omega)$ be a 
quadruple where $X$ is a quasi-projective variety with $G$ action, where $G$ is a reductive algebraic group, $L$ is a $G$-equivariant line bundle and $\omega$ is a $G$-invariant section of $L$. Our main example is $(\mathbb{A}^{n+1},\mathbb{C}^{\ast},\mathcal{O}(d),\omega)$. The $\mathbb{C}^{\ast}$ action on $\mathbb{A}^{n+1}$ is given by $\lambda\cdot(x_{0},x_{1},\cdots, x_{n})=(\lambda^{a_{0}}\cdot x_{0},\lambda^{a_{1}}\cdot x_{1},\cdots, \lambda^{a_{n}}\cdot x_{n})$, $a_{0}, a_{1}, \cdots, a_{n}$ are integers such that $\operatorname{gcd}(a_{0},a_{1}, \cdots, a_{n})=1$. $\mathcal{O}(d)$ is the trivial line bundle twisted with the character $\mathbb{X}_{d}: \mathbb{C}^{\ast}\rightarrow \mathbb{C}^{\ast}, \quad \lambda \mapsto \lambda^{d}$. 
  $\omega$ is a $\mathbb{C}^{\ast}$-invariant section of $\mathcal{O}(d)$. Namely $\omega$ is a degree $d$ polynomial, the weight of variable $x_{i}$ is $a_{i}$. 
  

We have $dg$ category $\operatorname{Fact}(X,G,L,\omega)$, whose objects are a quadruple $(\mathcal{E}_{-1},\mathcal{E}_{0},\Phi_{-1},\Phi_{0})$, where $\mathcal{E}_{-1}$ and $\mathcal{E}_{0}$ are $G$-equivariant quasi-coherent sheaves, $\Phi_{-1}:\mathcal{E}_{0}\rightarrow \mathcal{E}_{-1}\otimes L$ and $\Phi_{0}:\mathcal{E}_{-1}\rightarrow \mathcal{E}_{0}$ are morphism of $G$-equivariant sheaves such that
\begin{align*}
    \Phi_{-1}\circ \Phi_{0}=\omega.\\
    (\Phi_{0}\otimes L)\circ \Phi_{-1}=\omega.
\end{align*}
The space of morphisms in $\operatorname{Fact}(X,G,L,\omega)$ are the internal Hom of $G$-equivariant sheaves while extending the pairs of morphisms to certain $\mathbb{Z}$-graded complexes. We point out the reference \cite{BFK} for interested readers. There is a category $\operatorname{Acyclic}(\operatorname{Fact}(X,G,L,\omega))$ which imitates acyclic complexes in the category of complexes of sheaves. The absolute derived category $D^{abs}[\operatorname{Fact}(X,G,L,\omega)]$ is the homotopy category of $dg$ quotient $\frac{\operatorname{Fact}(X,G,L,\omega)}{\operatorname{Acyclic}(\operatorname{Fact}(X,G,L,\omega))}\in \operatorname{Hqe(dg-cat)}$. Let $\operatorname{Inj}(X,G,L,\omega)\subset \operatorname{Fact}(X,G,L,\omega)$ be the dg sub-category whose components are $G$-equivariant injective quasi-coherent sheaves. We write [$\mathcal{A}$] as the homotopic category of any dg category $\mathcal{A}$.
\begin{lem}
   The composition $\operatorname{Inj}(X,G,L,\omega)\rightarrow \operatorname{Fact}(X,G,L,\omega)\rightarrow \frac{\operatorname{Fact}(X,G,L,\omega)}{\operatorname{Acyclic}(\operatorname{Fact}(X,G,L,\omega))}$
   induces an equivalence of homotopic categories
   $$[\operatorname{Inj}(X,G,L,\omega)]\cong [\frac{\operatorname{Fact}(X,G,L,\omega)}{\operatorname{Acyclic}(\operatorname{Fact}(X,G,L,\omega))}]:=D^{abs}[\operatorname{Fact}(X,G,L,\omega)]$$
\end{lem}
Let $\operatorname{Inj_{coh}}(X,G,L,\omega)\subset \operatorname{Inj}(X,G,L,\omega)$ be a $dg$ sub-category whose objects are quasi-isomorphic to objects with coherent components in category $\operatorname{Fact}(X,G,L,\omega)$.
  
  Define shifting functor $$[1]: (\mathcal{E}_{-1},\mathcal{E}_{0},\Phi_{-1},\Phi_{0})\mapsto (\mathcal{E}_{0},\mathcal{E}_{-1}\otimes L,-\Phi_{0},-\Phi_{-1}\otimes L).$$ 
With cone construction, the homotopic category $[\operatorname{Inj_{coh}}(X,G,L,\omega)]$ is a triangulated category which is equivalent to the category of graded matrix factorization in \cite{orlov2009derived} for $(\mathbb{A}^{n+1},\mathbb{C}^{\ast},\mathcal{O}(d),\omega)$.
  
Denote by  
  $$\{1\}=-\otimes \mathcal{O}_{\mathbb{A}^{n+1}}(1):\operatorname{Inj_{coh}(\mathbb{A}^{n+1},\mathbb{C}^{\ast},\mathcal{O}(d),\omega)}\rightarrow \operatorname{Inj_{coh}(\mathbb{A}^{n+1},\mathbb{C}^{\ast},\mathcal{O}(d),\omega)}$$
  the twisting functor that maps 
$$\xymatrix{\mathcal{E}_{-1}
\ar[r]^{\Phi_{0}}&\mathcal{E}_{0}\ar[r]^{\Phi_{-1}}&\mathcal{E}_{-1}(d)}$$
to 
$$\xymatrix{\mathcal{E}_{-1}(1)
\ar[r]^{\Phi_{0}(1)}&\mathcal{E}_{0}(1)\ar[r]^{\Phi_{-1}(1)}&\mathcal{E}_{-1}(d+1)}$$
Clearly, we have equality of functors $\{d\}:=\{1\}^{d}=[2]$.
\par
 Let $X\subset \mathbb{P}(a_{1}, a_{2}, \cdots, a_{n})$ be a smooth hypersurface of degree $d\leq n$ defined by $\omega$. Let 
  $$\Ku(X):=\Big\langle \mathcal{O}_{X},\mathcal{O}_{X}(1),\cdots, \mathcal{O}_{X}(\sum^{n}_{j=0}a_{j}-1-d)\Big\rangle^{\perp}.$$
Roughly, $\Ku(X)$ is identified with the essential subcategory of $B$-branes of $X$ in Physics. If $X$ a is Calabi-Yau variety, $Ku(X)=D^{b}(X)$. On $\mathrm{LG}$ side, the category $[\operatorname{Inj_{coh}}(\mathbb{A}^{n+1},\mathbb{C}^{\ast},\mathcal{O}(d),\omega)]$ is identified with the category of $B$-branes of Landau-Ginzburg model. Physically, B-branes of $X$ and $LG$ model are naturally equivalent, which was proved by Orlov \cite{orlov2009derived} mathematically. Namely, we have equivalence 
$$\Ku(X)\cong[\operatorname{Inj_{coh}}(X,G,L,\omega)].$$
  Consider the natural enhancement $\operatorname{Inj_{coh}}(X)$, and let $\Ku_{dg}(X)$ be a $dg$ subcategory that enhance $\Ku(X)$. Orlov's $\sigma$ /LG correspondence can be lifted to be the equivalence of dg categories. 
\begin{thm}\label{prop_Orlov_identification}\cite[Theorem 6.13]{BFK}
  There is an equivalence in $\operatorname{Hqe(dg-cat)}$,
  $$\Phi:\operatorname{Inj_{coh}}(\mathbb{A}^{n+1},\mathbb{C}^{\ast},\omega)\cong \Ku_{dg}(X).$$
\end{thm}
According to \cite{BFK}, the natural functors can be reinterpreted as kernels of Fourier-Mukai transforms, and the natural transformations between these functors are morphism of kernels. We write $\Delta(m)$ as the kernel of functor $-\otimes\mathcal{O}_{\mathbb{A}^{n+1}}(m)$. 

\begin{lem}\label{serrefunctor}\cite[Theorem 1.2]{FeveroKellyLGfractionalcy}
The Serre functor of $[\operatorname{Inj_{coh}}(\mathbb{A}^{n+1},\mathbb{C}^{\ast},\omega)]$ is $-\otimes \mathcal{O}_{\mathbb{A}^{n+1}}(\sum^{n}_{j=0}-a_{j})[n+1]$.      
\end{lem}

\begin{proof}
 Since $\operatorname{Inj_{coh}}(\mathbb{A}^{n+1},\mathbb{C}^{\ast},\mathcal{O}(d), \omega)\cong \Ku_{dg}(X)$, the category $\operatorname{Inj_{coh}}(\mathbb{A}^{n+1},\mathbb{C}^{\ast}, \mathcal{O}(d),\omega)$ is smooth and proper, or by \cite[lemma 2.11, 2.14]{FeveroKellyLGfractionalcy}. Then there is a smooth proper dg algebra $A$ such that [$\operatorname{Inj_{coh}}(\mathbb{A}^{n+1},\mathbb{C}^{\ast},\mathcal{O}(d),\omega)]\cong D^{perf}(A)$, the arguments in \cite[Theorem 2.18]{FeveroKellyLGfractionalcy} show the Serre functor is $$-\otimes \omega_{\mathbb{A}^{n+1}}[n+1-\operatorname{dim}\mathbb{C}^{\ast}+1]=-\otimes \mathcal{O}_{\mathbb{A}^{n+1}}(\sum^{n}_{j=0}-a_{j})[n+1].$$
\end{proof}
Next we recall a key theorem in \cite[Theorem 1.2]{BFK}. For $g\in\mathbb{C}^{\ast}$, we write $W_{g}$ as the conormal sheaf of fixed locus $(\mathbb{A}^{n+1})^{g}$ and $k_{g}$ the character of $\operatorname{det}(W_{g})$. We write $H^{\bullet}(d\omega_{g})$ as the Koszul cohomology of the Jacobian ideal of $\omega_{g}:=\omega\vert_{(\mathbb{A}^{n+1})^{g}}$. Let $\gamma=e^{\frac{2\pi i}{d}}$, and $\mu_{d}=\langle1, \gamma, \gamma^{2}, \cdots, \gamma^{d-1}\rangle$. 
\begin{prop}\cite[Theorem 5.39]{BFK}\label{extendedhochschild}
 Assume $\omega$ has only isolated singularty at $0$, then 
\begin{align*}
 \Hom(\Delta,\Delta(m)[t])\cong& (\bigoplus_{g\in \mu_{d},\  t-\operatorname{rk}W_{g}\  \textbf{is even} }\mathrm{Jac}(\omega_{g})(m-k_{g}+d(\frac{t-\operatorname{rk}W_{g}}{2})))^{\mathbb{C}^{\ast}} 
 \end{align*}
\end{prop}
We refer the reader to \cite{BFK} for details of computation. We describe the multiplication under the isomorphism of Theorem~\ref{extendedhochschild}. To make this self contain, we introduce some notions used in the proof.
\par
Let $Z$ be a quasi-projective variety with $G$ action, $G$ is an algebraic group. Let $H$ be a closed subgroup of $G$. We have an action of $H$ on $G\times Z$ defined by
$$\tau: H\times G\times Z \rightarrow G\times Z,\quad (h,g,z)\mapsto (g\cdot h^{-1}, h\cdot z).$$ 
The fppf quotient of $G\times Z$ of $H$ is a scheme 
, which is denoted as $G\times ^{H}Z$, see \cite[Lemma 2.16]{BFK}. Consider morphisms
$$l: Z\rightarrow G\times^{H} Z,\quad x\mapsto (e,x).$$
$$\alpha: G\times^{H} Z\rightarrow Z, \quad (g,x)\mapsto gx.$$
First, the pull back functor $l^{\ast}$ define an equivalence of equivariant quasi-coherent sheaves. Namely 
$$l^{\ast}: \operatorname{Qcoh}_{G}G\times^{H}Z\rightarrow \operatorname{Qcoh}_{H}Z.$$
is an equivalence \cite[Lemma 1.3]{Tho97}. We write $\alpha_{\ast}: \operatorname{Qcoh}_{G}G\times^{H}Z\rightarrow \operatorname{Qcoh}_{G}Z$ as the push forward functor of $\alpha$, and $\alpha^{\ast}: \operatorname{Qcoh}_{G}Z \rightarrow\operatorname{Qcoh}_{G}G\times^{H}Z$ as the pull back functor of $\alpha$.
\begin{defn}
   $$\operatorname{Ind^{G}_{H}}:=\alpha_{\ast}\circ (l^{\ast})^{-1}: \operatorname{Qcoh}_{H}Z\rightarrow \operatorname{Qcoh}_{G}Z.$$
   $$\operatorname{Res}^{G}_{H}:=l^{\ast}\circ \alpha^{\ast}: \operatorname{Qcoh}_{G}Z\rightarrow \operatorname{Qcoh}_{H}Z.$$
\end{defn}
We still write $\operatorname{Ind}^{G}_{H}$ and $\operatorname{Res}^{G}_{H}$ as derived functors of derived categories of equivariant sheaves. $\operatorname{Ind}^{G}_{H}$ is right adjoint functor of $\operatorname{Res}^{G}_{H}$. In our case, $Z=\mathbb{A}^{n+1}\times \mathbb{A}^{n+1}$,  $G=\mathbb{C}^{\ast}\times_{\mathbb{C}^{\ast}}\mathbb{C}^{\ast}=\{(g_{1},g_{2})\in \mathbb{C}^{\ast}\times \mathbb{C}^{\ast}\vert g^{d}_{1}=g^{d}_{2}\}$, and $H=\mathbb{C}^{\ast}\subset G$ via diagonal embedding. The $G$ action on $Z$ is given by 
$$(g_{1},g_{2})\cdot (x_{1},\cdots, x_{n+1},y_{1},\cdots,y_{n+1})=(g^{a_{1}}_{1}x_{1},\cdots,g^{a_{n+1}}_{1}x_{n+1},g^{a_{1}}_{2}y_{1},\cdots,g^{a_{n+1}}_{2}y_{n+1})$$
By definition,
$$\Hom(\Delta, \Delta(m)[t])\cong 
\Hom_{[\operatorname{Inj_{coh}}(\mathbb{A}^{n+1}\times \mathbb{A}^{n+1}, \mathbb{C}^{\ast}\times_{\mathbb{C}^{\ast}} \mathbb{C}^{\ast}, -\omega\boxtimes \omega)]}(\operatorname{Ind}^{\mathbb{C}^{\ast}\times_{\mathbb{C}^{\ast}}\mathbb{C}^{\ast}}_{\mathbb{C}^{\ast}}\Delta_{\ast}\mathcal{O}_{\mathbb{A}^{n+1}},\operatorname{Ind}^{\mathbb{C}^{\ast}\times_{\mathbb{C}^{\ast}} \mathbb{C}^{\ast}}_{\mathbb{C}^{\ast}}\Delta_{\ast}\mathcal{O}_{\mathbb{A}^{n+1}}(m)[t]).$$
The multiplication
$$\Phi: \Hom(\Delta,\Delta(m_{1})[t_{1}])\times \Hom(\Delta,\Delta(m_{2})[t_{2}])\rightarrow \Hom(\Delta,\Delta(m_{1}+m_{2})[t_{1}+t_{2}])
$$
maps $(a,b)$ to $ab$ 
is the composition
\begin{equation}\xymatrix{\operatorname{Ind}_{\mathbb{C}^{\ast}}^{\mathbb{C}^{\ast}\times_{\mathbb{C}^{\ast}}\mathbb{C}^{\ast}}\Delta_{\ast}\mathcal{O}_{\mathbb{A}^{n+1}}\ar[r]^{b}&\operatorname{Ind}_{\mathbb{C}^{\ast}}^{\mathbb{C}^{\ast}\times_{\mathbb{C}^{\ast}}\mathbb{C}^{\ast}}\Delta_{\ast}\mathcal{O}_{\mathbb{A}^{n+1}}(m_{2})[t_{2}]\ar[r]^{a}&\operatorname{Ind}_{\mathbb{C}^{\ast}}^{\mathbb{C}^{\ast}\times_{\mathbb{C}^{\ast}}\mathbb{C}^{\ast}}\Delta_{\ast}\mathcal{O}_{\mathbb{A}^{n+1}}(m_{1}+m_{2})[t_{1}+t_{2}]}.
\end{equation}
To avoid cluttering the notation, we temporarily assume $m_{1}=m_{2}=t_{1}=t_{2}=0$. The sequence (2) is equivalent to
\begin{equation*}
\xymatrix@C=2cm{\mathbb{L}\Delta^{\ast}\operatorname{Res}^{\mathbb{C}^{\ast}\times_{\mathbb{C}^{\ast}}\mathbb{C}^{\ast}}_{\mathbb{C}^{\ast}}\operatorname{Ind}_{\mathbb{C}^{\ast}}^{\mathbb{C}^{\ast}\times_{\mathbb{C}^{\ast}}\mathbb{C}^{\ast}}\Delta_{\ast}\mathcal{O}_{\mathbb{A}^{n+1}}\ar[r]^{b}&\mathbb{L}\Delta^{\ast}\operatorname{Res}^{\mathbb{C}^{\ast}\times_{\mathbb{C}^{\ast}}\mathbb{C}^{\ast}}_{\mathbb{C}^{\ast}}\operatorname{Ind}_{\mathbb{C}^{\ast}}^{\mathbb{C}^{\ast}\times_{\mathbb{C}^{\ast}}\mathbb{C}^{\ast}}\Delta_{\ast}\mathcal{O}_{\mathbb{A}^{n+1}}\ar[r]^{a}&\mathcal{O}_{\mathbb{A}^{n+1}}}   
\end{equation*}
Here $a,b$ are regarded as $\mathbb{C}^{\ast}$ equivariant morphisms in $[\operatorname{Inj_{coh}}(\mathbb{A}^{n+1},\mathbb{C}^{\ast},0)]$ via diagonal embedding $\mathbb{C}^{\ast}\hookrightarrow \mathbb{C}^{\ast}\times_{\mathbb{C}^{\ast}} \mathbb{C}^{\ast}$.
The morphism $b$ here is $\mathbb{L}\Delta^{\ast}\operatorname{Res}^{\mathbb{C}^{\ast}\times_{\mathbb{C}^{\ast}}\mathbb{C}^{\ast}}_{\mathbb{C}^{\ast}}\Delta^{\ast}(b):\mathbb{L}\Delta^{\ast}\operatorname{Res}^{\mathbb{C}^{\ast}\times_{\mathbb{C}^{\ast}}\mathbb{C}^{\ast}}_{\mathbb{C}^{\ast}}\operatorname{Ind}_{\mathbb{C}^{\ast}}^{\mathbb{C}^{\ast}\times_{\mathbb{C}^{\ast}}\mathbb{C}^{\ast}}\mathcal{O}_{\mathbb{A}^{n+1}}\rightarrow \mathbb{L}\Delta^{\ast}\operatorname{Res}^{\mathbb{C}^{\ast}\times_{\mathbb{C}^{\ast}}\mathbb{C}^{\ast}}_{\mathbb{C}^{\ast}}\operatorname{Ind}_{\mathbb{C}^{\ast}}^{\mathbb{C}^{\ast}\times_{\mathbb{C}^{\ast}}\mathbb{C}^{\ast}}\Delta_{\ast}\mathcal{O}_{\mathbb{A}^{n+1}}$. 
\par

Next, let $\Gamma_{g}$ be the graph $x\mapsto (g\cdot x,x)$. The $\mathbb{C}^{\ast}\times_{\mathbb{C}^{\ast}}\mathbb{C}^{\ast}$ action on $\bigoplus_{g\in \mu_{d}}\mathcal{O}_{\Gamma_{g}}$ is defined by $(g_{1},g_{2})\cdot (gx,\cdot x_{1})=(g_{1}g\cdot x,g_{2}\cdot x)$. By \cite[Lemma 5.31]{BFK}, $\operatorname{Ind}_{\mathbb{C}^{\ast}}^{\mathbb{C}^{\ast}\times_{\mathbb{C}^{\ast}}\mathbb{C}^{\ast}}\Delta_{\ast}\mathcal{O}_{\mathbb{A}^{n+1}}\cong \bigoplus_{g\in \mu_{d}}\mathcal{O}_{\Gamma_{g}}$ as $\mathbb{C}^{\ast}\times_{\mathbb{C}^{\ast}}\mathbb{C}^{\ast}$ sheaves. Note that the decomposition of the object $\operatorname{Ind}_{\mathbb{C}^{\ast}}^{\mathbb{C}^{\ast}\times_{\mathbb{C}^{\ast}}\mathbb{C}^{\ast}}\Delta_{\ast}\mathcal{O}_{\mathbb{A}^{n+1}}$ into a direct sum corresponds exactly to the decomposition in Proposition~\ref{extendedhochschild}, which identifies different Hom-spaces with the pieces of Jacobian rings.  

\par
Let $a=(f_{1}, f_{\gamma}, \cdots, f_{\gamma^{j}}, \cdots f_{\gamma^{d-1}})\in \Hom(\Delta,\Delta(m_{1})[t_{1}])\cong \bigoplus^{d-1}_{j=0} \mathrm{Jac}(\omega_{\gamma^{j}})_{m_{1}-k_{\gamma^{j}}+d\frac{t_{1}-rk W_{\gamma^{j}}}{2}}$, and $b=~(g_{1}, g_{\gamma}, \cdots, g_{\gamma^{j}},\cdots, g_{\gamma^{d-1}})\in \Hom(\Delta,\Delta(m_{2})[t_{2}])\cong \bigoplus^{d-1}_{j=0}\mathrm{Jac}(\omega_{\gamma^{j}})_{m_{2}-k_{\gamma^{j}}+d\frac{t_{2}-rk W_{\gamma^{j}}}{2}}$. 
\begin{thm}\label{composition}
The multiplication map
$$\Phi: \Hom(\Delta,\Delta(m_{1})[t_{1}])\times \Hom(\Delta,\Delta(m_{2})[t_{2}])\rightarrow \Hom(\Delta,\Delta(m_{1}+m_{2})[t_{1}+t_{2}]), (a,b)\mapsto ab.$$

is given by the composition of the following diagram
\begin{equation}
 \xymatrix@C=3cm{\mathbb{L}\Delta^{\ast}\mathcal{O}_{\Gamma_{1}}\ar[r]^{g_{1}}\ar[dr]\ar[ddr]\ar[dddr]&\mathbb{L}\Delta^{\ast}\mathcal{O}_{\Gamma_{1}}(m_{2})[t_{2}]\ar[r]^{f_{1}}&\mathcal{O}_{\Gamma_{1}}(m_{1}+m_{2})[t_{1}+t_{2}]\\
\mathbb{L}\Delta^{\ast}\mathcal{O}_{\Gamma_{\gamma^{1}}}\ar[ur]^{g_{\gamma^{1}}}\ar[r]\ar[dr]\ar[ddr]&\mathbb{L}\Delta^{\ast}\mathcal{O}_{\Gamma_{\gamma^{1}}}(m_{2})[t_{2}]\ar[ur]^{f_{\gamma^{1}}}&\\
 \cdots& \cdots\ar[uur]&\\
\mathbb{L}\Delta^{\ast} \mathcal{O}_{\Gamma_{\gamma^{d-1}}}\ar[uuur]_{g_{\gamma^{d-1}}}\ar[r]&\mathbb{L}\Delta^{\ast}\mathcal{O}_{\Gamma_{\gamma^{d-1}}}(m_{2})[t_{2}]\ar[uuur]_{f_{\gamma^{d-1}}}&\\}   
\end{equation}
Here in the left box the component morphism $\mathbb{L}\Delta_{\ast}\mathcal{O}_{\Gamma_{\gamma_{i}}}\rightarrow \mathbb{L}\Delta_{\ast}\mathcal{O}_{\Gamma_{\gamma_{j}}}$ in $[\operatorname{Inj_{coh}}(\mathbb{A}^{n+1},\mathbb{C}^{\ast}, 0)]$ is the result of the action of the group
element $(\gamma^{j},1)$ on the morphism $g_{\gamma^{i-j}}$.
In particular $g_{1}\circ f_{1}\in \mathrm{Jac}(\omega)$ is multiplication of functions.

\end{thm}
\begin{proof} This is essentially the duality of functors $\operatorname{Res}^{\mathbb{C}^{\ast}\times_{\mathbb{C}^{\ast}}\mathbb{C}^{\ast}}_{\mathbb{C}^{\ast}}$ and $\operatorname{Ind^{\mathbb{C}^{\ast}\times_{\mathbb{C}^{\ast}}\mathbb{C}^{\ast}}_{\mathbb{C}^{\ast}}}$.
 The element $(\gamma^{k},1)\in \mathbb{C}^{\ast}\times_{\mathbb{C}^{\ast}} \mathbb{C}^{\ast}$ defines an isomorphism
$$\Hom(\mathcal{O}_{\Gamma_{\gamma^{i}}},\mathcal{O}_{\Gamma_{\gamma^{j}}})\cong \Hom(\mathcal{O}_{\Gamma_{\gamma^{i+k}}},\mathcal{O}_{\Gamma_{\gamma^{j+k}}}).$$
Since $b$ is a $\mathbb{C}^{\ast}\times_{\mathbb{C}^{\ast}}\mathbb{C}^{\ast}$ invariant morphism, other morphisms except $(g_{1},g_{\gamma},\cdots,g_{\gamma^{d-1}})$ in the left box are uniquely determined by the $\mathbb{C}^{\ast}$ invariant morphisms $g_{\gamma^{\bullet}}$ via diagonal embedding. After identifying $\Hom(\mathbb{L}\Delta^{\ast}\Delta_{\ast}\mathcal{O}_{\mathbb{A}^{n+1}}, \mathcal{O}_{\mathbb{A}^{n+1}})$ with certain homogeneous degree of $\mathrm{Jac}(\omega)$, $g_{1}\circ f_{1}$ is the composition of functions, hence multiplication of polynomials.
\end{proof}
\begin{rem}\label{compatibility} 
It is easy to observe that the Hochschild-Serre algebra of the graded matrix factorization is not commutative in general.
\end{rem}
\section{Kuznetsov's Fano threefold conjecture for quartic double solids and Gushel-Mukai threefolds}\label{section_main_theorem_proved}
\begin{thm}\label{thm_infinitesimal_torelli_for_Ku_fails}
Let $Y$ be a smooth quartic double solid, whose semi-orthogonal decomposition is given by 
$$D^b(Y)=\langle\Ku(Y),\oh_Y,\oh_Y(1)\rangle,$$
where $\Ku(Y)$ is the Kuznetsov component of the quartic double solid $Y$.  The canonical map $\gamma_{Y}$ induced by multiplication map (1) of Hochschild-Serre algebra 
$$\gamma_{Y}: \HH^{2}(\Ku(Y))\longrightarrow \Hom(\HH_{-1}(\Ku(Y)), \HH_{1}(\Ku(Y))).$$
has one dimensional kernel.
\end{thm}


    

\begin{proof}
We regard $Y$ as a degree $4$ smooth hypersurface in weighted projective space $\mathbb{P}(1,1,1,1,2)$. According to Theorem~\ref{prop_Orlov_identification}, $\Ku_{dg}(Y)\cong \operatorname{Inj_{coh}}(\mathbb{A}^{5},\mathbb{C}^{\ast}, \mathcal{O}(d),\omega)$, where $\omega$ is the polynomial defining $Y$, and the $\mathbb{C}^{\ast}$-action on $(x_0,x_1,x_2,x_3,x_4)$
is of weight $(1,1,1,1,2)$. Then by Proposition~\ref{extendedhochschild}, we have 
\begin{equation}
\Hom(\Delta,\Delta(m)[t])\cong (\bigoplus_{g\in\mu_{4},\  t-\operatorname{rk}W_{g}\  \textbf{is even} }\mathrm{Jac}(\omega_{g})(m-k_{g}+d(\frac{t-\operatorname{rk}W_{g}}{2})))^{\mathbb{C}^{\ast}},
\end{equation}
where $\mu_{4}=\{1,i,-1,-i\}$. 
\begin{itemize}
   \item If $g=1$, then $(\mathbb{A}^{5})^{g}=\mathbb{A}^{5}$, $\operatorname{rk}W_{g}=0, k_{g}=0$.
    \item If $g=i$, then $(\mathbb{A}^{5})^{g}=(0,0,0,0,0)$, $\operatorname{rk}W_{g}=5, k_{g}=-6$.
    \item If $g=-1$, then $(\mathbb{A}^{5})^{g}=(0,0,0,0,x_{5})$, $\operatorname{rk}W_{g}=4, k_{g}=-4$.
    \item If $g=-i$, then $(\mathbb{A}^{5})^{g}=(0.0,0,0,0)$, $\operatorname{rk}W_{g}=5, k_{g}=-6$.
\end{itemize}
Note that the Serre functor of the matrix factorization category is $-\otimes\mathcal{O}_{\mathbb{A}^5}(-6)[5]$ by Lemma~\ref{serrefunctor}, We write $\omega=x^{2}_{5}+f(x_{1}, x_{2}, x_{3}, x_{4})$, then
\begin{align}
    \HH_{-1}(\Ku(Y))\cong &\Hom(\Delta, \Delta(-6)[4])\cong \mathrm{Jac}(\omega)_{2}\oplus 0\oplus \mathrm{Jac}(\omega_{-1})_{-2}\oplus 0=\mathrm{Jac}(\omega)_{2}\\
    \HH_{1}(\Ku(Y))\cong &\Hom(\Delta, \Delta(-6)[6])\cong  \mathrm{Jac}(\omega)_{6}\oplus 0\oplus \mathrm{Jac}(\omega_{-1})_{2}\oplus 0=\mathrm{Jac}(\omega)_{6}\\
    \HH^{2}(\Ku(Y))\cong &\Hom(\Delta, \Delta[2])\cong \mathrm{Jac}(\omega)_{4}\oplus 0\oplus \mathrm{Jac}(\omega_{-1})_{0}\oplus 0=\mathrm{Jac}(\omega)_{4}\oplus k.
\end{align}
Let $(g_{1}, g_{i}, g_{-1}, g_{-i})\in \HH_{-1}(\Ku(Y))$ and $(f_{1}, f_{i}, f_{-1}, f_{-i})\in \HH^{2}(\Ku(Y))$, corresponding to (5) and (7) respectively. 
As explain in Theorem~\ref{composition} and see also the diagram (3), we have morphisms $(i,1)\cdot g_{-i}\in \Hom(\mathbb{L}\Delta^{\ast}\mathcal{O}_{\Gamma_{1}}, \mathbb{L}\Delta^{\ast}\mathcal{O}_{\Gamma_{i}})$, $(-1,1)\cdot g_{-1}\in \Hom(\mathbb{L}\Delta^{\ast}\mathcal{O}_{\Gamma_{1}},\mathbb{L}\Delta^{\ast}\mathcal{O}_{\Gamma_{-1}})$, $(-i,1)\cdot g_{i}\in \Hom(\mathbb{L}\Delta^{\ast}\mathcal{O}_{\Gamma_{1}}, \mathbb{L}\Delta^{\ast}\mathcal{O}_{\Gamma_{-i}})$; $(i,1)\cdot g_{1}\in \Hom(\mathbb{L}\Delta^{\ast}\mathcal{O}_{\Gamma_{i}},\mathbb{L}\Delta^{\ast}\mathcal{O}_{\Gamma_{i}})$, $(-1,1)\cdot g_{-i}\in \Hom(\mathbb{L}\Delta^{\ast}\mathcal{O}_{\Gamma_{i}},\mathbb{L}\Delta^{\ast}\mathcal{O}_{\Gamma_{-1}})$, $(-i,1)\cdot g_{-1}\in \Hom(\mathbb{L}\Delta^{\ast}\mathcal{O}_{\Gamma_{i}},\mathbb{L}\Delta^{\ast}\mathcal{O}_{\Gamma_{-i}})$; $(i,1)\cdot g_{i}\in \Hom(\mathbb{L}\Delta^{\ast}\mathcal{O}_{\Gamma_{-1}}, \mathbb{L}\Delta^{\ast}\mathcal{O}_{\Gamma_{i}})$, $(-1,1)\cdot g_{1}\in \Hom(\mathbb{L}\Delta^{\ast}\mathcal{O}_{\Gamma_{-1}},\mathbb{L}\Delta^{\ast}\mathcal{O}_{\Gamma_{-1}})$, $(-i,1)\cdot g_{-i}\in \Hom(\mathbb{L}\Delta^{\ast}\mathcal{O}_{\Gamma_{-1}}, \mathbb{L}\Delta^{\ast}\mathcal{O}_{\Gamma_{-i}})$; $(i,1)\cdot g_{-1}\in \Hom(\mathbb{L}\Delta^{\ast}\mathcal{O}_{\Gamma_{-i}}, \mathbb{L}\Delta^{\ast}\mathcal{O}_{\Gamma_{i}})$, $(-1,1)\cdot g_{i}\in \Hom(\mathbb{L}\Delta^{\ast}\mathcal{O}_{\Gamma_{-i}},\mathbb{L}\Delta^{\ast}\mathcal{O}_{\Gamma_{-1}})$, $(-i,1)\cdot g_{1}\in \Hom(\mathbb{L}\Delta^{\ast}\mathcal{O}_{\Gamma_{-i}}, \mathbb{L}\Delta^{\ast}\mathcal{O}_{\Gamma_{-i}})$.
According to Theorem~\ref{composition}, the composition 
$$\xymatrix{\HH^{2}(\Ku(Y))\times \HH_{-1}(\Ku(Y))\ar[r]&\HH_{1}(\Ku(Y))}.$$
is 
represented by
$$\xymatrix@C=4cm{1&\bullet\ar[r]^{g_{1}}\ar[dr]\ar[ddr]\ar[dddr]&\bullet\ar[r]^{f_{1}}&\bullet\\
i&\bullet\ar[r]\ar[dr]\ar[ddr]\ar[ur]^{g_{i}}&\bullet\ar[ur]^{f_{i}}&\\
-1&\bullet\ar[r]\ar[dr]\ar[ur]\ar[uur]^{g_{-1}}&\bullet\ar[uur]^{f_{-1}}&\\
-i&\bullet\ar[r]\ar[ur]\ar[uur]\ar[uuur]_{g_{-i}}&\bullet\ar[uuur]_{f_{-i}}&}.$$
Namely, the composition $(f_{1},f_{i},f_{-1},f_{-i})\circ(g_{1},g_{i},g_{-1},g_{-i})$ is 
\begin{align*}
(f_{1}\circ g_{1}+f_{i}\circ ((i,1)\cdot g_{-i}) +f_{-1}\circ((-1,1)\cdot g_{-1})+ f_{-i}\circ((-i,1)\cdot g_{i}),&\\
f_{1}\circ g_{i}+f_{i}\circ((i,1)\cdot g_{1})+ f_{-1}\circ ((-1,1)\cdot g_{-i})+ f_{-i}\circ  ((-i,1)\cdot g_{-1}),&\\
f_{1}\circ g_{-1}+f_{i}\circ ((i,1)\cdot g_{i}) +f_{-1}\circ ((-1,1)\cdot g_{1}) + f_{-i}\circ ((-i,1)\cdot g_{-i}) ,&\\
f_{1}\circ g_{-i}+ f_{i}\circ ((i,1)\cdot g_{-1})+f_{-1}\circ  ((-1,1)\cdot g_{i}) + f_{-i}\circ ((-i,1)\cdot g_{1}))
\end{align*}


Consider element $a=(0, 0, f_{-1}, 0)\in \HH^{2}(\Ku(Y))\cong \mathrm{Jac}(\omega)_{4}\oplus 0 \oplus \mathrm{Jac}(\omega_{-1})_{0}\oplus 0=\mathrm{Jac}(\omega)_{4}\oplus k$ and $b=(g_{1},0,0,0)\in \HH_{-1}(\Ku(Y))\cong \mathrm{Jac}(\omega)_{2}\oplus 0\oplus 0\oplus 0$. Then, 
$$(0, 0, f_{-1}, 0)\circ(g_{1}, 0, 0, 0)=(0\circ g_{1}+ f_{-1}\circ 0, 0, f_{-1}\circ ((-1,1)\cdot g_{1}), 0)=0\in \HH_{1}(\Ku(Y))\cong \mathrm{Jac}(\omega)_{6},$$
where the composition $f_{-1}\circ ((-1,1)\cdot g_{1})$ lies in $\mathrm{Jac}(\omega_{-1})_2=0$ by $(6)$, thus $f_{-1}\circ ((-1,1)\cdot g_{1})=0$. 

By \cite[Theorem 2.6]{Don} the map $$\mathrm{Jac}(f)_{4}\times \mathrm{Jac}(f)_{2}\rightarrow \mathrm{Jac}(f)_{6}.$$ is a non-degeneration multiplication. Thus the map $$\mathrm{Jac}(f)_{4}\rightarrow \Hom(\mathrm{Jac}(f)_{2},\mathrm{Jac}(f_6)).$$ is injective. On the other hand, simple computation shows $\mathrm{Jac}(\omega)_{4}=\mathrm{Jac}(f)_{4}$, $\mathrm{Jac}(\omega)_{2}=\mathrm{Jac}(f)_{2}$, and $\mathrm{Jac}(\omega)_{6}=\mathrm{Jac}(f)_{6}$. Then the map $$\mathrm{Jac}(\omega)_{4}\rightarrow \Hom(\mathrm{Jac}(\omega)_{2},\mathrm{Jac}(\omega)_6)$$ is also injective. 

Hence the canonical map 
$$\gamma_{Y}: \HH^{2}(\Ku(Y))\rightarrow \Hom(\HH_{-1}(\Ku(Y)),\HH_{1}(\Ku(Y))).$$
has one dimensional kernel.
\end{proof}
\begin{lem}\cite[Theorem 4.6]{jacovskis2022infinitesimal}\label{GMinfinitesimal}
    Let X be an ordinary GM threefold. Then the natural map 
    $$\gamma_{X}: \HH^{2}(\Ku(X))\rightarrow \Hom(\HH_{-1}(\Ku(X)),\HH_{1}(\Ku(X))).$$
    is injective.
\end{lem}
\begin{proof}
The map $\gamma$ in \cite[Theorem 4.6]{jacovskis2022infinitesimal} is related to $\gamma_{X}$ as 
$$\xymatrix@C=3cm{\HH^{2}(\Ku(X))\ar[r]^{\gamma_{X}}\ar[rd]^{\gamma}&\Hom(\HH_{-1}(\Ku(X)),\HH_{1}(\Ku(X))\ar[d]^{\simeq}\\
&\Hom(H^{2,1}(X),H^{1,2}(X))}.$$
Since $\gamma$ is injective, $\gamma_{X}$ is injective.
\end{proof}

\begin{lem}\label{lemma_injectivity_gamma}
Let $X$ be a special Gushel-Mukai threefold, then the morphism 
$$\gamma_X: \mathrm{HH}^2(\Ku(X))\rightarrow\mathrm{Hom}(\mathrm{HH}_{-1}(\Ku(X)),\mathrm{HH}_1(\Ku(X)))$$
is injective. 
\end{lem}

\begin{proof}
By \cite[Lemma 3.8]{kuznetsov2018derived} and \cite[Theorem 1.6]{kuznetsov2019categorical}, for any special Gushel-Mukai threefold $X$, there is an ordinary Gushel-Mukai threefold $X'$ such that $\Phi\colon\Ku(X')\simeq\Ku(X)$ is a Fourier-Mukai type equivalence. Then by \cite[Theorem 4.8]{jacovskis2022infinitesimal}, injectivity of $\gamma_X$ is equivalent to injectivity of $\gamma_{X'}$. By Lemma \ref{GMinfinitesimal}, $\gamma_{X'}$ is injective. Thus $\gamma_X$ is injective. 
\end{proof}

\begin{cor}\label{cor_main_theorem}
 For any Gushel-Mukai threefold $X$ and quartic double solid $Y$, there is no Fourier-Mukai type equivalence between the category $\Ku(X)$ and $\Ku(Y)$.
\end{cor}
\begin{proof}
Assume there is a Fourier-Mukai type equivalence $\Phi:\Ku(Y)\simeq\Ku(X)$ for any quartic double solid $Y$ and ordinary Gushel-Mukai threefold $X$. Then \cite[Theorem 4.8]{jacovskis2022infinitesimal} tells us the morphism $\gamma_{X}$ is injective if and only if $\gamma_Y$ is injective. Then by Lemma~\ref{GMinfinitesimal} and Lemma~\ref{lemma_injectivity_gamma}, the map $\gamma_X$ is injective for all smooth Gushel-Mukai threefolds. Thus $\gamma_Y$ is also injective, which contradicts Theorem~\ref{thm_infinitesimal_torelli_for_Ku_fails}. 
\end{proof}

\begin{rem}
In this paper, we work with dg-enhanced Kuznetsov categories, so any equivalence between them amounts to a Fourier-Mukai type equivalence. But in the cases of interest in this paper, all the equivalences between triangulated categories $\Ku(X)$ and $\Ku(Y)$ are proved to be of Fourier-Mukai type in \cite{li2022derived}, so there is no harm to work with enhanced Kuznetsov components.   
\end{rem}



\bibliographystyle{alpha}
{\small{\bibliography{reference}}}

\newcommand{\etalchar}[1]{$^{#1}$}
\begin{thebibliography}{CMHL{\etalchar{+}}23}

\bibitem[APR22]{altavilla2022moduli}
Matteo Altavilla, Marin Petkovic, and Franco Rota.
\newblock Moduli spaces on the kuznetsov component of fano threefolds of index 2.
\newblock {\em {\'E}pijournal de G{\'e}om{\'e}trie Alg{\'e}brique}, 6, 2022.

\bibitem[BFK14]{BFK}
Matthew Ballard, David Favero, and Ludmil Katzarkov.
\newblock A category of kernels for equivariant factorizations and its implications for hodge theory.
\newblock {\em Publ.math.IHES}, 120(1-111), 2014.

\bibitem[BFK23]{belmans2023hochschild}
Pieter Belmans, Lie Fu, and Andreas Krug.
\newblock Hochschild cohomology of hilbert schemes of points on surfaces.
\newblock {\em arXiv preprint arXiv:2309.06244}, 2023.

\bibitem[Bla16]{blanc_2016}
Anthony Blanc.
\newblock Topological k-theory of complex noncommutative spaces.
\newblock {\em Compositio Mathematica}, 152(3):489–555, 2016.

\bibitem[BLMS17]{bayer2017stability}
Arend Bayer, Mart{\'\i} Lahoz, Emanuele Macr{\`\i}, and Paolo Stellari.
\newblock {Stability conditions on Kuznetsov components}.
\newblock {\em (Appendix joint with Xiaolei Zhao) To appear in Ann. Sci. {\'E}c. Norm. Sup{\'e}r., arXiv:1703.10839}, 2017.

\bibitem[BMMS12]{bernardara2012categorical}
Marcello Bernardara, Emanuele Macr{\`\i}, Sukhendu Mehrotra, and Paolo Stellari.
\newblock A categorical invariant for cubic threefolds.
\newblock {\em Advances in Mathematics}, 229(2):770--803, 2012.

\bibitem[BO01]{bondal2001reconstruction}
Alexei Bondal and Dmitri Orlov.
\newblock Reconstruction of a variety from the derived category and groups of autoequivalences.
\newblock {\em Compositio Mathematica}, 125(3):327--344, 2001.

\bibitem[BP23]{bayer2022kuznetsov}
A.~Bayer and A.~Perry.
\newblock {Kuznetsov's {F}ano threefold conjecture via {K}3 categories and enhanced group actions}.
\newblock {\em J. Reine Angew. Math.}, 800:107--153, 2023.

\bibitem[Cal03]{caldararu2003mukai}
Andrei Caldararu.
\newblock The mukai pairing, i: the hochschild structure.
\newblock {\em arXiv preprint math/0308079}, 2003.

\bibitem[C{\u{a}}l05]{cualduararu2005derived}
Andrei C{\u{a}}ld{\u{a}}raru.
\newblock Derived categories of sheaves: a skimming.
\newblock {\em Snowbird lectures in algebraic geometry}, 388:43--75, 2005.

\bibitem[CMHL{\etalchar{+}}23]{CMHLZZ2024categorical}
Sebastian Casalaina-Martin, Xianyu Hu, Lin.Xun, Shizhuo Zhang, and Zheng Zhang.
\newblock Categorical torelli theorem via equivariant kuznetsov component of cubic threefolds.
\newblock {\em preprint}, 2023.

\bibitem[DJR23]{dell2023cyclic}
Hannah Dell, Augustinas Jacovskis, and Franco Rota.
\newblock Cyclic covers: Hodge theory and categorical torelli theorems.
\newblock {\em arXiv preprint arXiv:2310.13651}, 2023.

\bibitem[Don83]{Don}
Ron. Donagi.
\newblock Generic torelli for projective hypersurfaces.
\newblock {\em Compositio Mathematica}, 50(2-3):323--353, 1983.

\bibitem[FK18]{FeveroKellyLGfractionalcy}
David Favero and Tyler Kelly.
\newblock Fractional calabi–yau categories from landau–ginzburg models.
\newblock {\em Algebraic Geometry}, 5(5):596–649, 2018.

\bibitem[GLZ22]{guo2022conics}
Hanfei Guo, Zhiyu Liu, and Shizhuo Zhang.
\newblock Conics on gushel-mukai fourfolds, epw sextics and bridgeland moduli spaces.
\newblock {\em arXiv preprint arXiv:2203.05442}, 2022.

\bibitem[JLLZ21]{jacovskis2021categorical}
A.~Jacovskis, X.~Lin, Z.~Liu, and S.~Zhang.
\newblock {Categorical Torelli theorems for Gushel-Mukai threefolds}.
\newblock {\em arXiv preprint arXiv:2108.02946}, 2021.

\bibitem[JLLZ22]{jacovskis2022infinitesimal}
Augustinas Jacovskis, Xun Lin, Zhiyu Liu, and Shizhuo Zhang.
\newblock Infinitesimal categorical torelli theorems for fano threefolds.
\newblock {\em arXiv preprint arXiv:2203.08187}, 2022.

\bibitem[Kel06]{keller2006differential}
Bernhard Keller.
\newblock On differential graded categories.
\newblock {\em arXiv preprint math/0601185}, 2006.

\bibitem[KP18]{kuznetsov2018derived}
Alexander Kuznetsov and Alexander Perry.
\newblock Derived categories of gushel--mukai varieties.
\newblock {\em Compositio Mathematica}, 154(7):1362--1406, 2018.

\bibitem[KP19]{kuznetsov2019categorical}
Alexander Kuznetsov and Alexander Perry.
\newblock Categorical cones and quadratic homological projective duality.
\newblock {\em arXiv preprint arXiv:1902.09824}, 2019.

\bibitem[Kuz09]{kuznetsov2009derived}
Alexander~Gennad'evich Kuznetsov.
\newblock Derived categories of fano threefolds.
\newblock {\em Proceedings of the Steklov Institute of Mathematics}, 264:110--122, 2009.

\bibitem[LPZ22]{li2022derived}
Chunyi Li, Laura Pertusi, and Xiaolei Zhao.
\newblock Derived categories of hearts on kuznetsov components.
\newblock {\em arXiv preprint arXiv:2203.13864}, 2022.

\bibitem[LXZ24]{LXZ2024Hodge}
Xun Lin, Fei Xie, and Shizhuo Zhang.
\newblock Singular category: Hodge theory and (birational) categorical torelli theorem for nodal fano threefolds.
\newblock {\em preprint}, 2024.

\bibitem[LZ23]{lin2023serre}
Xun Lin and Shizhuo Zhang.
\newblock Serre algebra, matrix factorization and categorical torelli theorem for hypersurfaces.
\newblock {\em arXiv preprint arXiv:2310.09927}, 2023.

\bibitem[Orl03]{orlov2003derived}
Dmitri~Olegovich Orlov.
\newblock Derived categories of coherent sheaves and equivalences between them.
\newblock {\em Russian Mathematical Surveys}, 58(3):511, 2003.

\bibitem[Orl09]{orlov2009derived}
Dmitri Orlov.
\newblock Derived categories of coherent sheaves and triangulated categories of singularities.
\newblock {\em Algebra, Arithmetic, and Geometry: Volume II: In Honor of Yu. I. Manin}, pages 503--531, 2009.

\bibitem[Per22]{perry2020integral}
Alexander Perry.
\newblock {The integral Hodge conjecture for two-dimensional Calabi--Yau categories}.
\newblock {\em Compositio Mathematica}, 158(2):287--333, 2022.

\bibitem[PY22]{pertusi:some-remarks-fano-threefolds-index-two}
L.~Pertusi and S.~Yang.
\newblock {Some remarks on {F}ano three-folds of index two and stability conditions}.
\newblock {\em Int. Math. Res. Not. IMRN}, (17):12940--12983, 2022.

\bibitem[Tho87]{Tho97}
R.W Thomason.
\newblock Equivariant resolution, linearization, and hilbert's fourteenth problem over arbitrary base schemes.
\newblock {\em Advances in Mathematics, Volume 65, Issue 1}, 1987.

\bibitem[VI99]{shafa}
Yuri Vasilevich~Prokhorov V.A.~Iskovskikh.
\newblock {\em {Algebraic Geometry V: Fano Varieties}}.
\newblock Algebraic geometry. Springer, 1999.

\bibitem[Zha20]{Zhang2020Kuzconjecture}
Shizhuo Zhang.
\newblock {Bridgeland moduli spaces and Kuznetsov's Fano threefold conjecture}.
\newblock {\em arXiv preprint 2012.12193}, 2020.

\end{thebibliography}

\end{document}